\documentclass{elsart3-1}

\usepackage{amssymb}

\usepackage{graphicx}

\usepackage[english,francais]{babel}

\newtheorem{theorem}{Theorem}
\newtheorem{lemma}[theorem]{Lemma}
\newtheorem{e-proposition}[theorem]{Proposition}

\newtheorem{e-definition}[theorem]{Definition}
\newtheorem{remark}[theorem]{Remark}
\newtheorem{example}[theorem]{Example}
\newtheorem{assumption}[theorem]{Assumption}

\renewcommand{\theequation}{\arabic{equation}}
\def\og{\leavevmode\raise.3ex\hbox{$\scriptscriptstyle\langle\!\langle$~}}
\def\fg{\leavevmode\raise.3ex\hbox{~$\!\scriptscriptstyle\,\rangle\!\rangle$}}

\journal{the Acad\'emie des sciences}
\begin{document}
\centerline{}
\begin{frontmatter}


\selectlanguage{english}
\title{Exponential self-similar mixing and \\ loss of regularity for continuity equations}


\selectlanguage{english}
\author[authorlabel1]{Giovanni Alberti},
\ead{galberti1@dm.unipi.it}
\author[authorlabel2]{Gianluca Crippa},
\ead{gianluca.crippa@unibas.ch}
\author[authorlabel3]{Anna L.~Mazzucato}
\ead{alm24@psu.edu}

\address[authorlabel1]{Dipartimento di Matematica, 
Universit\`a di Pisa, largo Pontecorvo 5, 56127 Pisa, Italy}
\address[authorlabel2]{Departement Mathematik und Informatik,
Universit\"at Basel, Rheinsprung 21,
CH-4051 Basel, Switzerland}
\address[authorlabel3]{Department of Mathematics, Penn State University,
McAllister Building, University Park, PA 16802, USA}


\medskip
\begin{center}
{\small Received *****; accepted after revision +++++\\
Presented by �����}
\end{center}

\begin{abstract}
\selectlanguage{english}
We consider the mixing behaviour of the solutions of the continuity equation associated with a divergence-free velocity field. In this announcement we sketch two explicit examples of exponential decay of the mixing scale of the solution, in case of Sobolev velocity fields, thus showing the optimality of known lower bounds. We also describe how to use such examples to construct solutions to the continuity equation with Sobolev but non-Lipschitz velocity field exhibiting instantaneous loss of any fractional Sobolev regularity. {\it To cite this article: A.
Name1, A. Name2, C. R. Acad. Sci. Paris, Ser. I 340 (2005).}

\vskip 0.5\baselineskip

\selectlanguage{francais}
\noindent{\bf R\'esum\'e} \vskip 0.5\baselineskip \noindent
{\bf M\'elange auto-similaire exponentiel et perte de r\'egularit\'e pour l'\'equation de continuit\'e.}
Nous \'etudions le comportement de m\'elange des solutions de l'\'equation de continuit\'e associ\'e \`a un champ de v\'elocit\'e \`a divergence nulle. Dans cette annonce on esquisse deux exemples explicites de d\'ecalage exponentiel de l'\'echelle de m\'elange de la solution; dans le cas des champs de v\'elocit\'e Sobolev on d\'emontre donc l'optimalit\'e des estimations d'en bas connues. On d\'ecrit de m\^eme comment utiliser tels exemples pour construire des solutions de l'\'equation de continuit\'e aux champs de v\'elocit\'e Sobolev mais pas Lipschitzien : ces solutions perdent imm\'ediatement toute r\'egularit\'e Sobolev fractionnaire. {\it Pour citer cet article~: A. Name1, A. Name2, C. R. Acad. Sci.
Paris, Ser. I 340 (2005).}

\end{abstract}
\end{frontmatter}


\selectlanguage{english}
\section*{Introduction}

Consider a passive scalar $\rho$ which is advected by a time-dependent, divergence-free velocity field $u$ on the two-dimensional torus, i.e., a solution $\rho$ of the continuity equation with velocity field $u$:
\begin{equation}\label{e:continuity}
\left\{
\begin{array}{l} 
\partial_t \rho + {\rm div}\, (u \rho) = 0, \\
\rho(0,\cdot) = \bar \rho,
\end{array}
\right. \qquad \textrm{on ${\mathbb R}^+ \times {\mathbb T}^2$.}
\end{equation} 
We assume that the initial datum $\bar \rho$ satisfies $\int_{{\mathbb T}^2} \bar \rho = 0$ (this condition is, at least formally, preserved under the time evolution) and we are interested in the mixing behaviour of $\rho(t,\cdot)$ as $t$ tends to $+\infty$.

\medskip

In order to quantify the level of ``mixedness'', two different notions of mixing scale are available in the literature. The first one is based on homogeneous negative Sobolev norms (see for instance \cite{doering,mmp}), the most common one being the $\dot{H}^{-1}$ norm, which should be viewed as a characteristic length of the mixing in the system (here and in the following we use the dot as in $\dot{H}^{-1}$ to denote the {\em homogeneous} Sobolev norm):

\begin{e-definition} The functional mixing scale of $\rho(t,\cdot)$ is $\| \rho(t,\cdot) \|_{\dot{H}^{-1} (\mathbb{T}^2)}$.
\end{e-definition}

The second notion (see \cite{bressan}) is geometric and has been introduced for solutions with value $\pm 1$:

\begin{e-definition} Given $0<\kappa<1/2$, the geometric mixing scale of $\rho(t,\cdot)$ is the infimum $\varepsilon(t)$ of all $\varepsilon >0$ such that for every $x \in {\mathbb T}^2$ there holds
\begin{equation}\label{e:geom}
\kappa \leq \frac{{\mathcal L}^2 \big( \{ \rho(t,\cdot) = 1 \} \cap B_\varepsilon(x) \big)}{{\mathcal L}^2 (B_\varepsilon(x))} \leq 1 - \kappa \,,
\end{equation}
where ${\mathcal L}^2$ denotes the two-dimensional Lebesgue measure and $B_\varepsilon(x)$ is the ball of radius $\varepsilon$ centered at $x$.
\end{e-definition}
The parameter $\kappa$ is fixed and plays a minor role in the definition. Informally, in order for $\rho(t,\cdot)$ to have geometric mixing scale $\varepsilon(t)$, we require that every ball of radius $\varepsilon(t)$ contains a substantial portion of both level sets $\{ \rho(t,\cdot) = 1 \}$ and $\{ \rho(t,\cdot) = -1 \}$.

Although, strictly speaking, the two notions above are not equivalent (see the discussion in \cite{llnmd}), they are heuristically very related, and indeed most of the results in this area are available considering any of the two definitions.

\medskip

The mixing process is studied in the literature under energetic
 constraints on the velocity field, that is, assuming that the velocity field is bounded with respect to some spatial norm, uniformly in time. We now briefly review some of the related literature (most of the results hold indeed in any space dimension):
\begin{itemize}
\item[(a)] The velocity field $u$ is bounded in $\dot{W}^{s,p}(\mathbb{T}^2)$ uniformly in time for some $s<1$ and $1\leq p\leq \infty$ (the case $s=0$, $p=2$, often referred to as energy-constrained flow, is of particular interest in applications). In this case the Cauchy problem for the continuity equation~(\ref{e:continuity}) is not uniquely solvable in general (see \cite{abc2,abc1}). It is therefore possible to find a velocity field and a bounded solution which is non-zero at the initial time, but is identically zero at some later time. This fact means that it is possible to have perfect mixing in finite time, as already observed in \cite{doering} and established in \cite{llnmd} for $s=0$, building on an example from \cite{depauw,bressan}.
\item[(b)] The velocity field $u$ is bounded in $\dot{W}^{1,p}(\mathbb{T}^2)$ uniformly in time for some $1\leq p\leq \infty$ (the case $p=2$, often referred to as enstrophy-constrained flow, is of particular interest in applications). In this case the results in \cite{dpl} ensure uniqueness for the Cauchy problem~(\ref{e:continuity}), while for $p>1$ the quantitative Lipschitz estimates for regular Lagrangian flows in \cite{cdl} provide an exponential lower bound on the geometric mixing scale, $\varepsilon(t) \geq C \exp(-ct)$. The extension to the borderline case $p=1$ is still open (see, however, \cite{bc}). For the same class of velocity fields, an exponential lower bound for the functional mixing scale, $\| \rho(t,\cdot) \|_{\dot{H}^{-1}} \geq C \exp(-ct)$, has been proved in \cite{ikx,seis}.
\item[(c)] For velocity fields bounded in $\dot{BV}(\mathbb{T}^2)$ uniformly in time, the results in \cite{amb} ensure uniqueness for the Cauchy problem~(\ref{e:continuity}). It was recognized in \cite{bressan} that an exponential decay of the geometric mixing scale can indeed be attained for velocity fields bounded in $\dot{BV}(\mathbb{T}^2)$ uniformly in time, and actually the same example works also for the functional mixing scale.
\item[(d)] The velocity field $u$ is bounded in $\dot{W}^{s,p} (\mathbb{T}^2)$ uniformly in time for some $s>1$ (here the case of interest in applications is that of palenstrophy-constrained flows, that is, $s=2$ and $p=2$). In this case there are no better lower bounds for the decay of the (functional or geometric) mixing scale than the exponential one obtained for $s = 1$. The common belief, supported also by the numerical simulations in \cite{llnmd,ikx}, is that this bound is optimal.
\end{itemize}

\medskip

In this note we sketch two examples in which both the functional and the geometric mixing scales decay exponentially with a velocity field bounded in $\dot{W}^{1,p}(\mathbb{T}^2)$ uniformly in time (for every $1\leq p < \infty$ in the first example, and for every $1\leq p \leq \infty$ in the second example, thus including the Lipschitz case). These results show the sharpness of the lower bounds in \cite{cdl,seis,ikx} (point~(b) above), for the full range $1\leq p\leq\infty$. Our examples can be seen as a Sobolev (or even Lipschitz) variant of the example in point~(c). 

Moreover, we describe how such constructions can be employed to obtain counterexamples to the propagation of any fractional Sobolev regularity for solutions to the continuity equation~(\ref{e:continuity}) in $\mathbb{R}^d$, when the velocity field belongs to Sobolev classes that do not embed in the Lipschitz space.

\medskip

After the completion of the present work we were made aware of a related result obtained independently by Yao and Zlato\v{s} \cite{zlatos}, which provides examples of mixing of general initial data, with a rate that is optimal in the range $1 \leq p \leq \bar p$ for some $\bar p >2$.

\section{Scaling analysis in a self-similar construction}
\label{s:scaling}

A conceivable procedure to mix an initial datum is through a self-similar evolution. We fix $s \geq 0$ and $1 \leq p \leq \infty$ and assume the following:
\begin{assumption}\label{ass}{\rm
There exist a velocity field $u_0$ and a  solution $\rho_0$ to~(\ref{e:continuity}), both defined for $0 \leq t \leq 1$, such that: 
\begin{itemize}
\item[(i)] $u_0$ is bounded in $\dot{W}^{s,p} ({\mathbb T}^2)$ uniformly in time;
\item[(ii)] $\rho_0$ is bounded; 
\item[(iii)] there exists a constant $\lambda >0$, with $1/\lambda$ an integer, such that 
\[
\rho_0(1, x) = \rho_0 \left( 0, \frac{x}{\lambda} \right) \,.
\]
\end{itemize}}
\end{assumption}
We will give an explicit example of a $u_0$ and a $\rho_0$ satisfying this assumption (with $s=1$ and $1 \leq p < \infty$) in Section \ref{sec:2}.
Given a positive parameter $\tau$ to be determined later, we set for each integer $n=1,2,\ldots$  and for $t \in [0,\tau^n]$:
\[
u_n(t,x) = \frac{\lambda^n}{\tau^n} \; u_0 \left( \frac{t}{\tau^n} , \frac{x}{\lambda^n} \right) \,, 
\quad \rho_n(t,x) = \rho_0 \left( \frac{t}{\tau^n} , \frac{x}{\lambda^n} \right) \,.
\] 
It is easy to check that $\rho_n$ solves~(\ref{e:continuity}) with velocity field $u_n$ and that, because of Assumption~\ref{ass}(iii),
\begin{equation}\label{e:glue}
\rho_n (\tau^n , x) = \rho_{n+1} (0,x) \,.
\end{equation}

\begin{figure}
\begin{center}
  \includegraphics[scale=1]{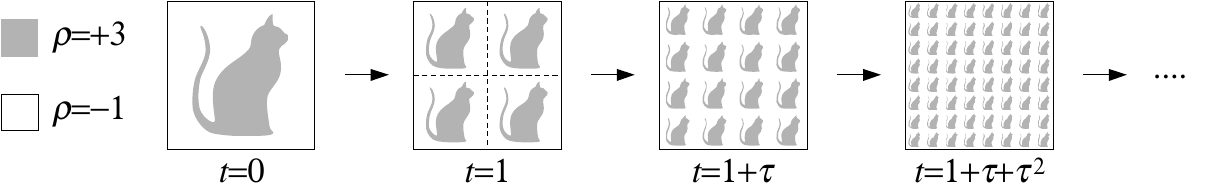}  
  \caption{Example of self-similar evolution for a function
$\rho$ taking only two values.} 
  \label{f:self}
\end{center}
\end{figure}

Now we construct $u$ and $\rho$ by patching together the velocity fields $u_0, u_1, \ldots$ and the corresponding solutions $\rho_0, \rho_1, \ldots$ In other words, we define
\[
u(t,x) = u_n (t - T_n , x) \,, \quad \rho(t,x) = \rho_n (t - T_n , x)  \quad \textrm{ for $T_n \leq t < T_{n+1}$, and $n=1,2,\ldots$,}
\]
where
\[
T_n = \sum_{i=0}^{n-1} \tau^i \,, \quad \textrm{ for $n=1, 2, \ldots, \infty$.}
\]
Observe that $u$ and $\rho$ are defined for $0 \leq t < T_\infty$, and that, thanks to equation~(\ref{e:glue}), $\rho$ is a solution of the Cauchy problem~(\ref{e:continuity}) 
with velocity field $u$ and initial condition $\bar{\rho}(x) = \rho_0(0,x)$. (See Figure \ref{f:self}.)

Since
\[
\left\| u_n (t,\cdot) \right\|_{\dot{W}^{s,p}(\mathbb{T}^2)} = \left( \frac{\lambda^{1-s}}{\tau} \right)^n \left\| u \left( \frac{t}{\tau^n} , \cdot \right) \right\|_{\dot{W}^{s,p}(\mathbb{T}^2)} \,,
\]
by setting $\tau = \lambda^{1-s}$ we obtain that $u$ is bounded in $\dot{W}^{s,p}({\mathbb T}^2)$ uniformly in time. Moreover, 
\[
\| \rho_n (t,\cdot) \|_{\dot{H}^{-1}(\mathbb{T}^2)} =  \lambda^n \left\| \rho_0 \left( \frac{t}{\tau^n} , \cdot \right) \right\|_{\dot{H}^{-1}(\mathbb{T}^2)} \leq M \lambda^n  \quad \textrm{ with } M = \sup_{0 \leq t \leq 1} \| \rho_0 (t,\cdot) \|_{\dot{H}^{-1}(\mathbb{T}^2)} \,,
\]
that is,
\begin{equation}\label{e:decay}
\| \rho(t,\cdot) \|_{\dot{H}^{-1}(\mathbb{T}^2)} \leq M \lambda^n, \qquad \textrm{ for $T_n \leq t < T_{n+1}$.}
\end{equation}

There are now three cases:
\begin{itemize}
\item[(1)] $s<1$, then $\tau = \lambda^{1-s} < 1$. In this case, $T_\infty$ is finite and $\| \rho(t,\cdot) \|_{\dot{H}^{-1}} \to 0$ as $t \to T_\infty$, that is, we have perfect mixing in finite time; 
\item[(2)] $s=1$, then $\tau = 1$. In this case, $T_\infty = \infty$ and $T_n = n$, therefore the inequality $t < T_{n+1}$ in~(\ref{e:decay}) becomes $t-1<n$, and the estimate in~(\ref{e:decay}) yields the following exponential decay of the functional mixing scale: 
\[
\| \rho(t,\cdot) \|_{\dot{H}^{-1}(\mathbb{T}^2)} \leq M \lambda^{t-1} \,;
\]
\item[(3)] $s>1$, then $\tau>1$. In this case, $T_\infty = \infty$ and
\[
T_n = \frac{\tau^n-1}{\tau - 1} = \frac{\lambda^{(1-s)n}-1}{\lambda^{1-s} - 1}.
\]
Hence, reasoning as above, (\ref{e:decay}) yields the following polynomial decay:
\[
\| \rho(t,\cdot) \|_{\dot{H}^{-1}(\mathbb{T}^2)} \leq M \frac{\left[ 1 + t (\lambda^{1-s} - 1) \right]^{-\frac{1}{s-1}}}{\lambda}
\simeq C(M,\lambda,s) \, t^{-\frac{1}{s-1}} \,.
\]
\end{itemize}

To summarize, we have proved the following:
 
\begin{theorem}\label{t:main}
Given $s$ and $p$, under Assumption~\ref{ass}, there exist a velocity field $u$ and a solution $\rho$ of the Cauchy problem (\ref{e:continuity}), such that $u$ is bounded in $\dot{W}^{s,p}(\mathbb{T}^2)$ uniformly in time and the functional mixing scale of $\rho$ exhibits the following behavior depending on $s$:
\begin{itemize}
\item[Case $s<1$:] perfect mixing in finite time;
\item[Case $s=1$:] exponential decay;
\item[Case $s>1$:] polynomial decay.
\end{itemize}
\end{theorem} 

\begin{remark}{\rm
In the above construction we can take $\rho_0$ with values $\pm 1$; in this case, using that for $T_n \leq t < T_{n-1}$ the solution $\rho(t,\cdot)$ is periodic on $\mathbb{T}^2$ with period $\lambda^n$, we can easily prove that the geometric mixing scale $\varepsilon(t)$ exhibits the same behavior of the functional mixing scale detailed in Theorem~\ref{t:main}. Note that for $s>1$ our self-similar examples do not match the known exponential lower bound for the (geometric and functional) mixing scale, which is supposed to be optimal.}
\end{remark}

\section{Two constructions} \label{sec:2}

It is possible to construct an example of a velocity field $u_0$ and a solution $\rho_0$ satisfying Assumption~\ref{ass} by considering a shear flow that splits a ball into two pieces, afterward rearranged into two smaller balls. This is essentially the nature of the example in \cite{bressan}. Note that such a shear flow is discontinuous along a line and indeed it is of class $BV$ but not of class $W^{1,1}$. 

 In order to obtain a Sobolev velocity field one could consider a ``modulated'' shear flow of the type $u(t,x_1,x_2) = (0,v(t,x_1))$, for which the ``total advection'' $x_1 \mapsto \int_0^1 v(t,x_1) \, dt$ is a characteristic function. However, such a construction cannot be achieved (even) with  $v \in L^1_t (\dot{H}^1_{x_1})$. And indeed the examples with Sobolev regularity
we present will be genuinely two-dimensional.

\medskip

The first example shows that Assumption~\ref{ass} can be fulfilled for $s=1$ and $1 \leq p < \infty$. We will actually only describe the time evolution of a set in $\mathbb{T}^2$ for $0 \leq t \leq 1$; then $\rho_0$ will be the function taking the value $+1$ on the set and $-1$ on the complement, and $u_0$ will be a velocity field that realizes such an evolution. Given the evolution of the set, the velocity field $u_0$ is provided by the following lemma, except at the time(s) when a singularity occurs in the geometric evolution:
\begin{lemma}\label{sel}
Given a smooth set evolving smoothly in time, such that the area of each connected component is constant in time, there exists a smooth, divergence-free velocity field the flow of which deforms the set according to the given evolution.
\end{lemma}
\begin{example}\label{ex:bottle}{\rm
We consider the geometric evolution described in Figure~\ref{f:pinching}; the velocity field $u_0$ realizing this evolution can be taken bounded in $\dot{W}^{1,p}(\mathbb{T}^2)$, for any given $1\leq p < \infty$, uniformly in time.
\begin{figure}[h]
\begin{center}
  \includegraphics[scale=1]{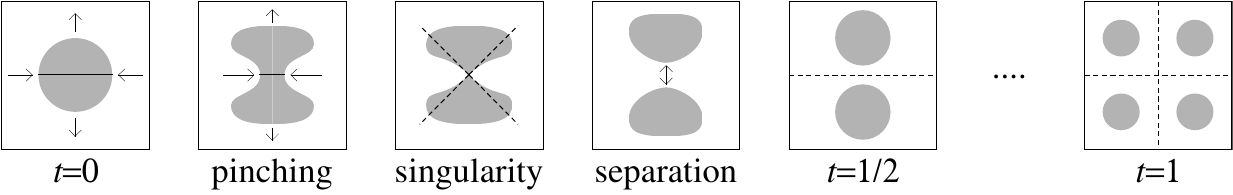}  
  \caption{First construction: splitting a connected component
in two.} 
  \label{f:pinching}
\end{center}
\end{figure}
The evolution is smooth except at the times at which a connected component of the set is split in two. This splitting is realized by first creating a small tube-like region in the connected component, which is then pinched, thus separating the component into two parts. One way to achieve the pinching with the correct regularity for the velocity field is to impose a horizontal compression with a square root-like velocity along the horizontal diameter of the ball, which is therefore shrunk to a single point in finite time. In order to obtain a (globally) divergence-free velocity field, this horizontal compression has to be compensated by a vertical decompression. The details of the construction are actually quite involved, therefore we leave them for a forthcoming paper.}
\end{example}

\begin{remark}{\rm
Although the results in \cite{dpl,amb} guarantee the uniqueness of a suitable ``flow solution'' for Sobolev velocity fields and for velocity fields with bounded variation, we see in Example~\ref{ex:bottle} that the set of initial data lacking {\em pointwise} uniqueness may contain a full line, and that the flow does not preserve topological properties of sets such as connectedness.}
\end{remark}


\medskip

Since the evolution in Example~\ref{ex:bottle} does not preserve connectedness, clearly it cannot be implemented with a Lipschitz velocity field. Our second example has Lipschitz regularity and strictly speaking is not self-similar; rather, it consists of two basic elements, which (up to rotations) are replicated at each step at finer scales. What we do is reminiscent of the construction of the classical Peano curve filling the square
(again, we leave the details for a forthcoming paper).

\begin{example}\label{ex:snake}{\rm
The basic moves of the construction are illustrated in Figure~\ref{f:snake}. We start from a straight strip at time $t=0$. Using Lemma~\ref{sel}, the strip is rearranged into the meandering strip at $t=1$, which consists of ``straight'' and ``curved'' pieces. 
\begin{figure}[h]
\begin{center}
  \includegraphics[scale=1]{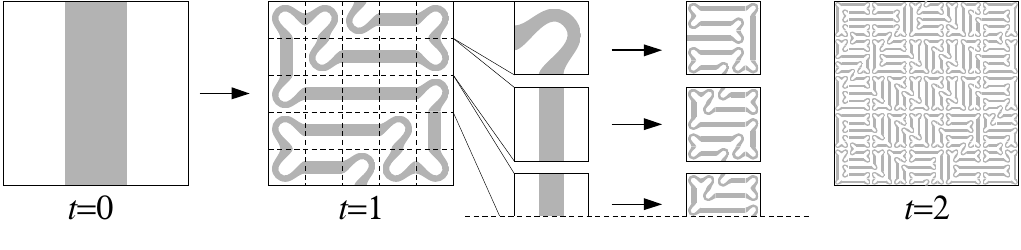}  
  \caption{Second construction: the Peano snake.} 
  \label{f:snake}
\end{center}
\end{figure}
In the time interval from $t=1$ to $t=2$ the set evolves as follows: the ``straight pieces'' are evolved employing the same transformation (possibly together with a rotation) used from $t=0$ to $t=1$, while the ``curved pieces'' are rearranged as shown in Figure~\ref{f:snake} for the right upper corner at $t=1$, again using Lemma~\ref{sel}. Notice that at every step our picture contains only ``straight pieces'' and ``curved pieces'', therefore the entire evolution is determined by the two basic moves in Figure~\ref{f:snake}.}
\end{example}

\begin{remark}{\rm
Although lacking exact self-similarity, the fact that only two basic elements are repeated allows to follow the analysis in Section~\ref{s:scaling} and to prove that both the functional and geometric mixing scales decay exponentially in time.}
\end{remark}

\section{Loss of regularity for solutions to the continuity equation}

We finally describe how to employ the above constructions to obtain counterexamples to the propagation of any fractional Sobolev regularity for solutions to the continuity equation~(\ref{e:continuity}) in $\mathbb{R}^d$. 

First of all, we observe that the construction in Example~\ref{ex:snake} can be modified to obtain a divergence-free velocity field $u \in C^\infty( {\mathbb R}^+ \times {\mathbb R}^d)$ and a 
solution $\rho \in C^\infty( {\mathbb R}^+ \times {\mathbb R}^d)$ to problem~(\ref{e:continuity}), so that
\begin{itemize}
\item[(i)] both $u$ and $\rho$ are supported in the unit cube $Q$
(in space);
\item[(ii)] $u$ is  bounded in $W^{1, \infty}(Q)$ uniformly in time;
\item[(iii)] $\| \rho(t,\cdot) \|_{\dot{H}^{-s}(Q)} \lesssim \exp(-c_s t)$
for all $s>0$ and some positive constant $c$.
\end{itemize}
Since the $L^2$ norm of $\rho$ is conserved in time, the interpolation inequality
\[
\| \rho(t,\cdot) \|_{L^2(Q)} \leq \| \rho(t,\cdot) \|_{\dot{H}^s(Q)}^{1/2} \| \rho(t,\cdot) \|_{\dot{H}^{-s}(Q)}^{1/2}
\]
gives
\begin{equation}\label{e:saturation}
\| \rho(t,\cdot) \|_{\dot{H}^{s}(Q)} \gtrsim \exp(c_s t) \quad \textrm{for all $s>0$.}
\end{equation}
This our example ``saturates'' the Gronwall estimate $\| \rho(t,\cdot) \|_{\dot{H}^{s}(Q)} \lesssim \exp(c_s t)$.

Then, given a disjoint family $Q_n$ of cubes in ${\mathbb R}^d$ with side $\lambda_n$, we can consider in each $Q_n$ (suitable translations of) the rescaled velocity fields and the solutions
\[
v_n (t,x) = \frac{\lambda_n}{\tau_n} \; u \left( \frac{t}{\tau_n} , \frac{x}{\lambda_n} \right)\,,
 \qquad
\theta_n (t,x) = C_n \; \rho\left( \frac{t}{\tau_n} , \frac{x}{\lambda_n} \right)\,,
\]
and define $v = \sum_n v_n$ and $\theta = \sum_n \theta_n$. It is immediate to check that such $\theta$ solves the continuity equation~(\ref{e:continuity}) with velocity field $v$. Moreover, thanks to the exponential growth in~(\ref{e:saturation}), we can
choose the construction parameters $\lambda_n, \tau_n, C_n$
so to impose a certain regularity of the velocity field $v$ and
of the initial datum $\theta(0,\cdot) = \bar\theta$ on
the one hand, and ensure an instantaneous loss of any
Sobolev regularity for the solution $\theta(t,\cdot)$ on the
other hand.
More precisely:
\begin{theorem}
We can choose the construction parameters so that:
\begin{itemize}
\item[(i)] the velocity field $v$ is bounded in $\dot W^{1,p}({\mathbb R}^d)$ uniformly in time, for any given $p<\infty$;
\item[(ii)] the initial datum $\bar\theta$ belongs to $C_c^\infty({\mathbb R}^d)$;
\item[(iii)] the solution $\theta(t,\cdot)$ does not belong to
$\dot{H}^s({\mathbb R}^d)$ for any $s>0$ and $t>0$.
\end{itemize}
Moreover both $v$ and $\theta$ are compactly supported
in space, and can be taken smooth on the complement of a point in ${\mathbb R}^d$.
\end{theorem}
Also in the case of velocity fields bounded in $\dot W^{r,p}({\mathbb R}^d)$ uniformly in time our construction gives some loss of regularity of the initial datum, provided that $\dot W^{r,p}({\mathbb R}^d)$ does not embed in $\dot W^{1,\infty}({\mathbb R}^d)$.



\section*{Acknowledgements}
This work was started during a visit of the first and third authors
at the University of Basel, whose kind hospitality is acknowledged. Their stay has been partially supported by the Swiss National Science Foundation grant 140232. The third author was partially supported by the US National Science Foundation grants  DMS 1009713, 1009714, 1312727.



\end{document}